\documentclass[11pt, reqno]{amsart}
\allowdisplaybreaks[3]

\usepackage[indentafter]{titlesec}
\usepackage{enumerate}
\usepackage{mathrsfs}
\usepackage{amsopn,amscd,amsmath,amssymb,amsthm,amsfonts,epsfig,graphics}
\usepackage{pslatex,avant,color}
\usepackage{cite}

\theoremstyle{plain}

\newtheorem{propo}{Proposition}[section]
\newtheorem{defi}{Definition}[section]
\newtheorem{exm}{Example}[section]

\theoremstyle{remark}
\newtheorem{rem}{Remark}[section]

\titleformat{name=\section}{}{\thetitle.}{0.8em}{\centering\scshape}
\titleformat{name=\subsection}{}{\thetitle.}{0.5em}{\bfseries}[]

\numberwithin{equation}{section}
\begin{document}

\title[Characterization of Jordan Vectors of Operator-Valued Functions]{Characterization of Jordan Vectors of Operator-Valued Functions with Applications in Differential Equations}
\author{Muhamed Borogovac}

\begin{abstract}
A well-known characterization of Jordan vectors of a matrix polynomial $L(z)$ is generalized to a characterization of Jordan vectors  of the operator-valued function $Q(z)$ at an eigenvalue $\alpha \in \mathbb{C}$. The results are then applied to solve a system of nonlinear ordinary differential equations.
\end{abstract}


\maketitle
\thispagestyle{empty}


\section{Introduction}\label{s2}

\textbf{1.1.} Let $\mathbb{N}$, $\mathbb{R}$, and $\mathbb{C}$ denote the sets of positive integers, real numbers, and complex numbers, respectively. The extended complex plane is defined as $\bar{\mathbb{C}}:=\mathbb{C}\cup \left\lbrace \infty \right \rbrace $. Let $\mathcal{H}$ denote a Hilbert space, and let $\mathcal{L}(\mathcal{H})$ denote the Banach space of bounded linear operators on $\mathcal{H}$. 

We assume that $\mathcal{D}(Q) \subseteq \mathbb{C}$ is a domain, and that $Q:\mathcal{D}(Q) \to \mathcal{L}(\mathcal{H})$ is a holomorphic operator-valued function on $\mathcal{D}(Q)$, that is, the domain of $Q$ and its domain of holomorphy coincide. We have chosen this level of generality so that our main result, Proposition \ref{proposition44}, applies to important classes of operator-valued functions, such as generalized Nevanlinna functions, which are meromorphic on $\mathbb{C} \setminus \mathbb{R}$ and may have generalized poles as well as continuous and residual spectra on the real axis; see \cite{Lu3}. We frequently focus on $\mathcal{H}= \mathbb{C}^{n}$, where $Q$ is a matrix-valued function and where we are able to concretely apply the abstract result of Proposition \ref{proposition44}.

If an operator-valued function $Q(z)$ is boundedly invertible at at least one point of holomorphy $w \in \mathcal{D}(Q)$, then it is boundedly invertible on some domain in $\mathbb{C}$ and is called a \textit{regular} operator-valued function on that domain; see, e.g., \cite[p.~327]{Lu4}. Evidently, a matrix-valued function holomorphic on some domain is regular if and only if $\chi(z) := \det Q(z) \not\equiv 0$.

Recall that an isolated singularity $\beta \in \mathbb{C}$ is called a \textit{pole of order} $m \in \mathbb{N}$ of the operator-valued function $Q(z)$ if $m$ is the minimal number such that the function $\tilde{Q}(z) := (z-\beta)^{m} Q(z)$ is holomorphic at $\beta$. We frequently deal with rational matrix-valued functions $Q$ that are \textit{meromorphic} on $\bar{\mathbb{C}}$, that is, they have only points of holomorphy and poles in $\bar{\mathbb{C}}$.
\\[2ex]
\textbf{1.2.} In this subsection, we review the terminology and basic facts related to matrix polynomials, since the concepts of eigenvalues, eigenvectors, and Jordan vectors introduced for matrix polynomials, as well as Proposition \ref{proposition14} below, will be extended to more general classes of operator-valued functions in this paper.

Let $A_{j},\, \, j=0,\, 1,\, \mathellipsis ,\, l$, be a set of $n\times n$ constant complex matrices, where $l \in \mathbb{N}$, and let 
\begin{equation}
\label{eq16}
A_{l}\frac{d^{l}\textbf{u}}{{dt}^{l}}+\mathellipsis +A_{1}\frac{d\textbf{u}}{dt}+A_{0}\textbf{u}=0
\end{equation}
be the corresponding linear system of differential equations, where 
\[
\textbf{u}(t):=\left(
\begin{array}{*{20}c}
u_{1}(t)\\
\vdots \\ 
u_{n}(t) \\ 
\end{array} \right)
\]
is an n-dimensional unknown complex vector function. This system is called a \textit{homogeneous linear system of ordinary differential equations (ODE) with constant coefficients}. One can seek a solution in the form
\begin{equation}
\label{eq110}
\textbf{u}\left( z \right)=\boldsymbol{\varphi} e^{zt}, \quad \boldsymbol{\varphi} :=\left( {\begin{array}{*{20}c}
\varphi_{1}\\
\vdots \\
\varphi_{n}\\
\end{array} } \right)\ne 0, \quad \varphi_{i}\in \mathbb{C},\quad i=1,2,\, \mathellipsis ,\, n.
\end{equation}
After substituting (\ref{eq110}) into (\ref{eq16}), we obtain:
\begin{equation}
\label{eq112}
L\left( z \right)\boldsymbol{\varphi} e^{zt}=0,
\end{equation}
where 
\[
L\left( z \right):=A_{l}z^{l}+\mathellipsis 
+A_{1}z+A_{0} .
\]
A vector $\boldsymbol{\varphi} \ne 0$ is a solution of (\ref{eq112}) for some $z=\alpha$, or equivalently, a solution of the algebraic system:
\begin{equation}
\label{eq114}
L\left( \alpha \right)\boldsymbol{\varphi}=0,
\end{equation}
if and only if $z=\alpha$ is a solutions of the equation:
\begin{equation}
\label{eq116}
\det {L\left( z \right)}=0.
\end{equation}
The polynomial $\chi (z) :=\det {L\left( z \right)}$ is called the \textit{characteristic polynomial} of the matrix polynomial $L(z)$. The zeros of the characteristic polynomial, i.e., solutions of the equation (\ref{eq116}) are called \textit{eigenvalues} of $L\left( z \right)$. Every solution $\boldsymbol{\varphi} \ne 0$ of (\ref{eq114}) is called an \textit{eigenvector} corresponding to the eigenvalue $\alpha $. 

The following proposition, see \cite[Proposition 1.9]{GLR}, is very important in the theory of ordinary linear differential equations.

\begin{propo}\label{proposition14} The vector function 
\begin{equation}
\label{eq118}
\textbf{u}\left( t \right)=\left( \frac{t^{k-1}}{(k-1)!}\boldsymbol{\varphi} 
_{0}+\frac{t^{k-2}}{(k-2)!}\boldsymbol{\varphi}_{1}+\mathellipsis +\boldsymbol{\varphi}_{k-1} 
\right)e^{\alpha t},
\end{equation}
where $k\in \mathbb{N}$ and  $\boldsymbol{\varphi}_{j}\in \mathbb{C}^{n}$ for $j=0,\, 1,\, \mathellipsis ,\, k-1$, is a solution of equation (\ref{eq16}) if and only if the following equalities hold: 
\begin{equation}
\label{eq120}
\sum\limits_{p=0}^i {\frac{1}{p!}L^{\left( p \right)}\left( \alpha 
\right)\boldsymbol{\varphi}_{i-p}} =0, \, i=0,\, \mathellipsis ,\, k-1,
\end{equation}
where $L^{\left( p \right)}\left( \alpha \right)$ is the p-th derivative of the matrix polynomial $L\left( z \right)$ at the eigenvalue $\alpha $ of $L\left( z \right)$.
\end{propo}
The sequence of n-dimensional vectors $\boldsymbol{\varphi}_{0},\, \boldsymbol{\varphi}_{1},\, \mathellipsis ,\, \boldsymbol{\varphi}_{m-1}$, with $\boldsymbol{\varphi}_{0}\ne 0$ and $m\leq k $, for which identities (\ref{eq120}) hold, is called a \textit{Jordan chain of order (or length) $m$} \textit{for the matrix polynomial} $L\left( z \right)$ \textit{at the eigenvalue} $\alpha \in \mathbb{C}$. The Jordan chain is said to be \textit{maximal, of order $k$}, if the system (\ref{eq120}) with the additional equation for $i=k$ does not have a solution. 
\\[2ex]
\textbf{1.3.} Evidently, rational matrix-valued functions, which are meromorphic on $\bar{\mathbb{C}}$, are more general than matrix polynomials, as they may also have poles in $\mathbb{C}$. Unlike matrix polynomials, for which all eigenvalues (including their multiplicities) are characterized as the solutions of the characteristic equation, the situation for a rational matrix-valued function $Q(z)$ is more complex: in addition to zeros, there are also poles of the characteristic function $\chi (z) := \det Q(z)$, and the zeros and poles of the characteristic function are not sufficient to characterize the zeros and poles of the matrix-valued function $Q$ itself or their multiplicities. For the purpose of characterizing the zeros and poles of operator-valued functions, including matrix-valued functions, the so-called root functions and pole cancellation functions are used.

The following definition is consistent with those given in \cite{GS,B1}. See also \cite[Definition 2.3]{GT}.

\begin{defi}\label{definition42}
Let $Q: \mathcal{D}(Q) \rightarrow \mathcal{L}(\mathcal{H})$ be an operator-valued function holomorphic on $\mathcal{D}(Q)$. Suppose that the vector function $\boldsymbol{\varphi}(z)$ satisfies the following conditions at some point $\alpha \in \mathbb{C}$:
\begin{enumerate}[(a)]
\item $\boldsymbol{\varphi}(\alpha) \ne 0$;
\item $(Q(z)\boldsymbol{\varphi}(z))^{(l)} \to 0$ as $z \to \alpha$, for $0 \le l \le m-1$, where $l, m-1 \in \mathbb{N} \cup \{0\}$.
\end{enumerate}
Then $\boldsymbol{\varphi}(z)$ is called a \textit{root function} of $Q(z)$ of order at least $m$ at the critical point (or zero)~$\alpha$. If $m$ is the maximal number for which condition~(b) holds, then $\boldsymbol{\varphi}(z)$ is said to be a root function of \textit{exact order} (or \textit{exact partial multiplicity}) $m$ at~$\alpha$.  

The vector $\boldsymbol{\varphi}(\alpha)$ is called an \textit{eigenvector corresponding to the eigenvalue}~$\alpha$. If $\boldsymbol{\varphi}(z)$ is a root function of maximal order $k$ among all root functions with the same eigenvector $\boldsymbol{\varphi}(\alpha)$, then $\boldsymbol{\varphi}(z)$ is called a \textit{canonical root function} associated with $\boldsymbol{\varphi}(\alpha)$.

The maximal order among all canonical root functions at~$\alpha$ is called the \textit{order of the zero} at~$\alpha$.
\end{defi}

If an operator-valued function $Q(z)$ is boundedly invertible in some domain, i.e., regular, then $\alpha$ is a pole of $Q(z)$ if and only if $\alpha$ is a zero of $Q^{-1}(z):=Q(z)^{-1}$. Thus, all the above definitions for a zero $\alpha$ of $Q$ can be applied to a pole $\alpha$ of $Q^{-1}$.

Our main goal is to find the necessary and sufficient conditions under which the system (\ref{eq120}) holds for the eigenvalues of a class of operator-valued functions $Q$ that are more general than matrix polynomials. In this general situation, we will not be able to use solutions of any system of differential equations the way the solutions (\ref{eq118}) were used in Proposition \ref{proposition14}. However, we will see that certain rational matrix-valued functions are associated with specific systems of nonlinear differential equations.

\section{Jordan vectors of an operator-valued function $Q\left( z \right)$}\label{s4}
\begin{propo}\label{proposition44}
Let $Q: \mathcal{D}(Q) \rightarrow \mathcal{L}(\mathcal{H})$ be an operator-valued function holomorphic on $\mathcal{D}(Q)$, and let $\alpha \in \mathcal{D}(Q)$ be a zero of $Q(z)$ of order at least $l \in \mathbb{N}$ (hence, $\alpha$ is not a singularity of $Q$).

A function $\boldsymbol{\varphi}(z)$, holomorphic at $\alpha$, satisfying $\boldsymbol{\varphi}(\alpha) \neq 0$, is a root function of 
$Q(z)$ at $\alpha$ of order at least $l \in \mathbb{N}$ if and only if it can be written in the form
\begin{equation}\label{eq42}
\boldsymbol{\varphi}(z) = \sum_{s=0}^{l-1} \frac{1}{s!}(z - \alpha)^{s} \boldsymbol{\varphi}^{(s)}(\alpha) + (z - \alpha)^{l-1} O(z- \alpha),
\end{equation}
where the derivatives
\begin{equation}\label{eq44}
Q^{(j)}(z), \,\, \boldsymbol{\varphi}^{(j)}(z), \,\, \forall j = 0, 1, \dots, l-1,
\end{equation}
satisfy the conditions
\begin{equation}\label{eq46}
\sum\limits_{p=0}^j \frac{1}{p!(j-p)!} Q^{(p)}(\alpha)\boldsymbol{\varphi}^{(j-p)}(\alpha) = 0, \quad j = 0, \mathellipsis, l-1.
\end{equation}
Here, $O(z- \alpha)$ is a function such that, for some $M>0$, it satisfies
\[
\mid O(z- \alpha)\mid \leq M \mid z-\alpha \mid, \quad as \quad z \rightarrow \alpha.
\]
\end{propo}
\begin {proof}
Assume that $\alpha \in \mathcal{D}(Q)$ is a zero of $Q(z)$ of order at least $l \in \mathbb{N}$, and that $\boldsymbol{\varphi}(z)$ is a root function of $Q$ at $\alpha$ of order at least $l$. Then, in some neighborhood of $\alpha$, the derivatives (\ref{eq44}) and the derivatives $(Q(z)\boldsymbol{\varphi}(z))^{(j)}(z)$, $\forall j = 0, 1, \ldots, l-1$, exist.

The Taylor expansion of the function $\psi(z) := Q(z)\boldsymbol{\varphi}(z)$ is
\[
Q(z)\boldsymbol{\varphi}(z) = \sum_{j=0}^{l-1} \frac{(z-\alpha)^j}{j!}\big(Q(z)\boldsymbol{\varphi}(z)\big)^{(j)}_{\mid z=\alpha} 
    + \sum_{j=l}^{\infty} \frac{(z-\alpha)^j}{j!}\big(Q(z)\boldsymbol{\varphi}(z)\big)^{(j)}_{\mid z=\alpha} 
\]
\[
= \sum_{j=0}^{l-1} \frac{(z-\alpha)^j}{j!}\big(Q(z)\boldsymbol{\varphi}(z)\big)^{(j)}_{\mid z=\alpha} 
    +(z-\alpha)^{l-1} \sum_{j=1}^{\infty} \frac{(z-\alpha)^j}{(j+l-1)!}\big(Q(z)\boldsymbol{\varphi}(z)\big)^{(j+l-1)}_{\mid z=\alpha}. 
\]
When we introduce the function 
\begin{equation}
\label{eq48}
O(z-\alpha):=\sum_{j=1}^{\infty} \frac{(z-\alpha)^j}{(j+l-1)!}\big(Q(z)\boldsymbol{\varphi}(z)\big)^{(j+l-1)}_{\mid z=\alpha} 
\end{equation}
and develop the derivatives 
\[
\big(Q(z)\boldsymbol{\varphi}(z)\big)^{(j)}_{\mid z=\alpha}=\sum\limits_{p=0}^j \frac{j!}{p!(j-p)!} Q^{(p)}(\alpha)\boldsymbol{\varphi}^{(j-p)}(\alpha) 
\]
we obtain
\begin{equation}
\label{eq49}
Q(z)\boldsymbol{\varphi}(z) = \sum_{j=0}^{l-1} \frac{(z-\alpha)^j}{j!} \sum\limits_{p=0}^j \frac{j!}{p!(j-p)!} Q^{(p)}(\alpha)\boldsymbol{\varphi}^{(j-p)}(\alpha) 
    + (z-\alpha)^{l-1} O(z-\alpha).
\end{equation}
According to Definition \ref{definition42} (b) satisfies
\[  
Q(z)\boldsymbol{\varphi}(z)= (z-\alpha)^{l-1} O(z-\alpha).
\]
This implies that the first sum in (\ref{eq49}) is equal to zero. Since the functions $(z-\alpha)^j$, $j = 0, \ldots, l-1$, are linearly independent, the condition (\ref{eq46}) is satisfied.

Conversely, assume that the condition (\ref{eq46}) is satisfied, i.e.,
\begin{equation}
\label{eq410}
\sum\limits_{p=0}^j \frac{j!}{p!(j-p)!} Q^{(p)}(\alpha)\boldsymbol{\varphi}^{(j-p)}(\alpha) = 0, \quad j = 0, \mathellipsis, l-1.
\end{equation}

By substituting the Taylor expansions of $Q(z)$ and $\boldsymbol{\varphi}(z)$ at $\alpha$ of degree $l$ into the function $\psi(z) := Q(z)\boldsymbol{\varphi}(z)$, and introducing the $O(z-\alpha)$ functions related to $Q(z)$ and $\varphi (z)$ analogously to (\ref{eq48}), we obtain
\[
Q(z)\boldsymbol{\varphi}(z) =
\]
\[
\Big(\sum_{p=0}^{l-1} \frac{(z-\alpha)^p}{p!} Q^{(p)}(\alpha) + (z-\alpha)^{l-1} O(z-\alpha)\Big)\times
\]  
\[  
  \times \Big(\sum_{s=0}^{l-1} \frac{(z-\alpha)^s}{s!} \boldsymbol{\varphi}^{(s)}(\alpha) + (z-\alpha)^{l-1} O(z-\alpha)\Big)
\]
\[
= \sum\limits_{j=0}^{l-1} \frac{1}{j!} (z-\alpha)^j 
   \sum\limits_{p=0}^j \frac{j!}{p!(j-p)!} Q^{(p)}(\alpha)\boldsymbol{\varphi}^{(j-p)}(\alpha)
   + (z-\alpha)^{l-1} O(z-\alpha).
\]
The last equality follows by introducing the exponent $j := p + s$ of $(z-\alpha)$ and grouping the products $Q^{(p)}(\alpha)\boldsymbol{\varphi}^{(j-p)}(\alpha)$ according to $(z-\alpha)^{j}$ for $j = 0, \ldots, l-1$, while collecting all terms with $j > l-1$ into $(z-\alpha)^{l-1}O(z-\alpha)$.

Since condition (\ref{eq410}) is satisfied, it follows that
\[
Q(z)\boldsymbol{\varphi}(z) = (z-\alpha)^{l-1} O(z-\alpha).
\]
This proves that the function $\boldsymbol{\varphi}(z)$ given by (\ref{eq42}) is indeed a root function of $Q$ at $\alpha$ of order at least $l$.\end{proof} 
\begin{defi}\label{definition46}
The vectors $0 \neq \boldsymbol{\varphi}_{j} \in \mathcal{H}$ that satisfy the conditions
\begin{equation}\label{eq412}
\sum\limits_{p=0}^j \frac{1}{p!} Q^{(p)}(\alpha) \boldsymbol{\varphi}_{j-p} = 0, \quad j = 0, \mathellipsis, l-1,
\end{equation}
are called generalized Jordan vectors of the operator-valued function $Q(z)$ at $\alpha$.

If $l \in \mathbb{N}$ is the maximal number for which the system (\ref{eq412}) has a solution, then the vectors $\boldsymbol{\varphi}_s$, $s = j - p = 0, \ldots, l - 1$, form the maximal Jordan chain corresponding to the eigenvector $\boldsymbol{\varphi}_0$.
\end{defi}

\begin{rem}\label{remark47}
Assume that $L(z)$ is a matrix polynomial. Then, evidently, conditions (\ref{eq46}) are satisfied by the Jordan vectors $\boldsymbol{\varphi}_0, \ldots, \boldsymbol{\varphi}_{l-1}$ of $L(z)$ at $\alpha$, that satisfy condition (\ref{eq120}). In other words, Proposition \ref{proposition44} is a \textbf{generalization of Proposition \ref{proposition14}}, extending the result from matrix polynomials to operator-valued functions including matrix-valued functions.

In this case, the root function $\boldsymbol{\varphi}(z)$ of $L(z)$ at $\alpha$ is given by (\ref{eq42}), with 
\[
\boldsymbol{\varphi}^{(s)}(\alpha) = s! \boldsymbol{\varphi}_s, \quad s = 0, \mathellipsis, l-1.
\]
\end{rem} \hfill $\square$
\begin{rem}\label{remark48}
If $A \in \mathcal{L}(\mathcal{H})$ and $\alpha \in \mathbb{C}$ is an eigenvalue of $A$, then substituting $Q(z) = A - zI$ into condition (\ref{eq46}) yields the well-known relations between the Jordan vectors $\boldsymbol{\varphi}_s$ of the operator $A$:
\begin{equation}\label{eq414}
(A - \alpha I)\boldsymbol{\varphi}_0 = 0, \quad \boldsymbol{\varphi}_{s-1} = (A - \alpha I)\boldsymbol{\varphi}_s, \quad s = 1, \ldots, l-1.
\end{equation}

If $\boldsymbol{\varphi}_0, \ldots, \boldsymbol{\varphi}_{l-1}$ is a Jordan chain, not necessarily maximal, of the operator $A \in \mathcal{L}(\mathcal{H})$ at the eigenvalue $\alpha \in \mathbb{C}$, then a root function of the operator-valued function $Q(z)=A - \alpha I$ of order $l$ is given by 
\[
\boldsymbol{\varphi} (z)= \boldsymbol{\varphi}_{0}+(z-\alpha)\boldsymbol{\varphi}_{1}+ ... +(z-\alpha)^{l-1}\boldsymbol{\varphi}_{l-1}.
\]
Indeed, 
\[
Q(z)\boldsymbol{\varphi} (z)=(A-zI)(\boldsymbol{\varphi}_{0}+(z-\alpha)\boldsymbol{\varphi}_{1}+ ... +(z-\alpha)^{l-1}\boldsymbol{\varphi}_{l-1})=
\]
\[
=\left[ (A-\alpha I)+(\alpha -z)I\right] (\boldsymbol{\varphi}_{0}+(z-\alpha)\boldsymbol{\varphi}_{1}+ ... +(z-\alpha)^{l-1}\boldsymbol{\varphi}_{l-1})= 
\]
\[
=(A-\alpha I)(\boldsymbol{\varphi}_{0}+(z-\alpha)\boldsymbol{\varphi}_{1}+ ... +(z-\alpha)^{l-1}\boldsymbol{\varphi}_{l-1})+
\]
\[
+(\alpha -z)(\boldsymbol{\varphi}_{0}+(z-\alpha)\boldsymbol{\varphi}_{1}+ ... +(z-\alpha)^{l-1}\boldsymbol{\varphi}_{l-1})=-(z-\alpha)^{l}\boldsymbol{\varphi}_{l-1}.
\]
The last equality follows after applying the relations (\ref{eq414}) to the first sum and subsequent cancellations with the corresponding members in the second sum.

Note that similar things are discussed in \cite[Subsection 2.2.]{GT} for a matrix $A$. \hfill $\square$
\end{rem}

Unlike scalar functions and matrix polynomials, the rational matrix-valued functions may have zeros and poles at the same point. In this case we cannot apply equations (\ref{eq46}) to find Jordan chains and the root functions at such points. In the following example we exhibit a such function $Q(z)$.
\begin{exm}
\label{example410} Find a root function at the singularity $\alpha =0$ of the matrix-valued function
\[
Q(z)=\left( {\begin{array}{*{20}c}
z &  2 & -1\\
0 &  z-1 & -2\\
0 &  0 & \frac{3}{z}\\
\end{array} } \right).
\]
\end{exm}
It is easy to verify that 
\[
\boldsymbol{\varphi}_{1}(z) = 
\begin{pmatrix}
1 \\[4pt]
z \\[4pt]
0
\end{pmatrix},\quad \boldsymbol{\varphi}_{2}(z) = 
\begin{pmatrix}
1 \\[4pt]
0 \\[4pt]
z^{2}
\end{pmatrix}, \quad \boldsymbol{\varphi}_{3}(z) = 
\begin{pmatrix}
1 \\[4pt]
0 \\[4pt]
0
\end{pmatrix},
\]
all linear combinations of those, and many other functions, are root functions of  $Q(z)$ at $\alpha = 0$. Hence, the singularity $\alpha = 0$ is both a zero and a pole of $Q$. \hfill $\square$

In the following example, we demonstrate how the equations (\ref{eq46}) can be applied to find generalized Jordan vectors and the corresponding root functions of the rational matrix-valued functions, i.e. that Proposition \ref{proposition44} is indeed an generalization of Proposition \ref{proposition14}.
\begin{exm}
\label{example411}
For the meromorphic matrix-valued function
\[
Q(z) = 
\begin{pmatrix}
\frac{z-2}{z-3} & \frac{1}{3-z} & 0 \\[4pt]
0 & \frac{z-2}{z-3} & \frac{1}{z-3} \\[4pt]
0 & 0 & \frac{z+3}{z-3}
\end{pmatrix},
\]
find the zeros $\alpha_1$ and $\alpha_2$, and the corresponding maximal root functions $\boldsymbol{\varphi}_1(z)$ and $\boldsymbol{\varphi}_2(z)$.
\end{exm}

Since, in this example, $\chi(z)$ and $Q(z)$ have identical zeros and poles, including their multiplicities, there are no common poles and zeros of $Q$. Therefore, the zeros of $Q$ are points of holomorphy. Hence, the conditions of Proposition \ref{proposition44} are satisfied for zeros of $Q$. It follows that:
\[
\det Q(z)=\frac{(z-2)^{2}(z+3)}{(z-3)^{3}}=0 \Rightarrow \big(\left( \alpha_{1}=2, k_{1}=2\right) \wedge \left( \alpha_{2}=-3, k_{2}=1\right)\big).
\] 
Let us find the generalized Jordan vectors $\boldsymbol{\varphi}(z)$ at $\alpha_{1} = 2$ of $Q$. Since $\alpha_{1} = 2$ is a zero of order $k_{1} = 2$, we will use the first two equations of the system (\ref{eq46}) to find the first two Jordan vectors: 
\begin{equation}
\label{eq416}
\begin{array}{*{20}c}
Q(2)\boldsymbol{\varphi}_{0} = 0,\\
Q'(2)\boldsymbol{\varphi}_{0} + Q(2)\boldsymbol{\varphi}_{1} = 0.\\
\end{array}
\end{equation}
We have
\[
Q(2) = 
\begin{pmatrix}
0 & 1 & 0\\
0 & 0 & -1\\
0 & 0 & -5
\end{pmatrix},
\]
\[
Q'(z) = 
\begin{pmatrix}
\frac{-1}{(z-3)^{2}} & \frac{1}{(z-3)^{2}} & 0\\
0 & \frac{-1}{(z-3)^{2}} & \frac{-1}{(z-3)^{2}}\\
0 & 0 & \frac{-6}{(z-3)^{2}}
\end{pmatrix}
\Rightarrow 
Q'(2) = 
\begin{pmatrix}
-1 & 1 & 0\\
0 & -1 & -1\\
0 & 0 & -6
\end{pmatrix}.
\]
The first equation of the system (\ref{eq416}) gives
\[
\begin{pmatrix}
0 & 1 & 0\\
0 & 0 & -1\\
0 & 0 & -5
\end{pmatrix}
\boldsymbol{\varphi}_{0} =
\begin{pmatrix}
0\\
0\\
0
\end{pmatrix}
\Rightarrow 
\boldsymbol{\varphi}_{0} =
\begin{pmatrix}
1\\
0\\
0
\end{pmatrix}.
\]
The second equation implies
\[
\begin{pmatrix}
-1 & 1 & 0\\
0 & -1 & -1\\
0 & 0 & -6
\end{pmatrix}
\begin{pmatrix}
1\\
0\\
0
\end{pmatrix}
+
\begin{pmatrix}
0 & 1 & 0\\
0 & 0 & -1\\
0 & 0 & -5
\end{pmatrix}
\boldsymbol{\varphi}_{1}
= 0. \quad We \quad select \quad
\boldsymbol{\varphi}_{1} =
\begin{pmatrix}
0\\
1\\
0
\end{pmatrix}.
\]
It is straightforward to verify that the third equation of the system does not have a solution consistent with the first two equations. Therefore, there are only two generalized eigenvectors. The root function at $\alpha_{1} = 2$ is 
\[
\boldsymbol{\varphi}(z) =
\begin{pmatrix}
1\\
0\\
0
\end{pmatrix}
+
\begin{pmatrix}
0\\
1\\
0
\end{pmatrix}
(z - 2).
\]
It is easy to check that
\[
Q(z)\boldsymbol{\varphi}(z) =
\begin{pmatrix}
0\\
\frac{(z - 2)^{2}}{z - 3}\\
0
\end{pmatrix},
\]
which confirms that $\alpha_{1} = 2$ is indeed a zero of $Q(z)$ of order two.

Finding the second root function, which is a root function of order $k_{2} = 1$, is easier. We need only the first equation of the system (\ref{eq46}). This simple task is left to the reader. \hfill $\square$

In the following example, we will see that the first equation of (\ref{eq46}), when $Q(z)$ is a rational matrix-valued function, can also be used to obtain a solution of a certain nonlinear system of differential equations.
\begin{exm}\label{example412} Solve the system of differential equations
\begin{equation}
\label{eq418} 
\begin{array}{*{20}c}
\frac{2}{u_{1}^{(2)}}+\frac{1}{u_{1}}+\frac{2}{u_{2}^{(2)}}+\frac{1}{u_{2}^{(1)}}=0,\\
\frac{-1}{u_{1}^{(2)}}+\frac{1}{u_{1}^{(1)}}+\frac{1}{u_{2}^{(2)}}=0.\\
\end{array}
\end{equation}
\end{exm}
We seek a solution of the form
\begin{equation}
\label{eq420}
\left( {\begin{array}{*{20}c}
u_{1}(t)\\
u_{2}(t)\\
\end{array} }\right) =\left( {\begin{array}{*{20}c}
\frac{1}{\varphi_{1}}\\
\frac{1}{\varphi_{2}}\\
\end{array} }\right)e^{z t}.
\end{equation}
By substituting (\ref{eq420}) into (\ref{eq418}), we obtain 
\begin{equation}
\label{eq422}
\left( {\begin{array}{*{20}c}
\frac{2}{z^{2}}+1 & \frac{2}{z^{2}}+\frac{1}{z} \\
\frac{-1}{z^{2}}+\frac{1}{z} & \frac{1}{z^{2}} \\
\end{array} }\right) \left( {\begin{array}{*{20}c}
\varphi_{1}\\
\varphi_{2}\\
\end{array} }\right)e^{-z t}=\left( {\begin{array}{*{20}c}
0\\
0\\
\end{array} }\right).
\end{equation} 
We have
\begin{equation}
\label{eq424}
Q(z)=\left( {\begin{array}{*{20}c}
\frac{2}{z^{2}}+1 & \frac{2}{z^{2}}+\frac{1}{z} \\
\frac{-1}{z^{2}}+\frac{1}{z} & \frac{1}{z^{2}} \\
\end{array} }\right) \, \wedge \, \det Q(z)= \frac{-z+4}{z^{4}}.
\end{equation}
We see that $Q(z)$ has a zero of order $k=1$ at $\alpha = 4$.
This means that the algebraic system (\ref{eq46}) reduces to its first equation.
After multiplying by $4^{2}$ to avoid fractions, we have
\[
Q(4)\boldsymbol{\varphi} =0 \Rightarrow 
\left( {\begin{array}{*{20}c}
18 & 6 \\
3 & 1 \\
\end{array} }\right) \left( {\begin{array}{*{20}c}
\varphi_{1}\\
\varphi_{2}\\
\end{array} }\right)=\left( {\begin{array}{*{20}c}
0\\
0\\
\end{array} }\right) \Rightarrow \varphi_{2}=-3\varphi_{1}.
\]
Hence, a solution (\ref{eq420}) of the system (\ref{eq418}) is given by
\[
\left( {\begin{array}{*{20}c}
u_{1}(t)\\
u_{2}(t)\\
\end{array} }\right) =\left( {\begin{array}{*{20}c}
\frac{1}{\varphi_{1}}\\
\frac{1}{-3\varphi_{1}}\\
\end{array} }\right)e^{4t}, \quad \varphi_{1}\in \mathbb{C}\setminus\lbrace 0 \rbrace. 
\] \hfill $\square$

In the previous example, we showed how to associate a system of nonlinear differential equations of the form (\ref{eq418}) with the rational matrix-valued function $Q(z)$ by seeking a solution of the form (\ref{eq420}). Then the solutions (\ref{eq420}) exist if and only if the vector $\boldsymbol{\varphi}$ is an eigenvector of the function $Q(z)$, i.e., a solution of the first equation of (\ref{eq46}). In fact, the following proposition holds.

\begin{rem}\label{remark414} Let $Q(z)$ be an $n\times n$ rational matrix-valued function associated with a system of the form (\ref{eq418}). This means that $n$ differential equations of order $l \in \mathbb{N}$ with $n$ unknown functions $u_{m}(t), m=1, ..., n$, are sums of fractions of the form $\frac{a_{ij}^{k}}{u^{(k)}_{m}}$, where $a_{ij}^{k}$, $ k=0,...,l;  \, i,j=1,...,n$, are complex constants and $u^{(k)}_{m}$ are derivatives of the unknown functions. The function $Q(z)$ is associated with the system by the substitutions 
\begin{equation}
\label{eq426}
\textbf{u}(t):=\left(
\begin{array}{*{20}c}
u_{1}(t)\\
\vdots \\ 
u_{n}(t) \\ 
\end{array} \right)=\left( {\begin{array}{*{20}c}
\frac{1}{\varphi_{1}}\\
\vdots \\ 
\frac{1}{\varphi_{n}}\\
\end{array} }\right)e^{z t}
\end{equation}
as demonstrated in the previous example. Let $z=\alpha$ be an eigenvalue of $Q(z)$. Then the function (\ref{eq426}), for $z=\alpha$, is a solution of the system if and only if the vector
\[
\boldsymbol{\varphi}=\left(
\begin{array}{*{20}c}
\varphi_{1}\\
\vdots \\ 
\varphi_{n} \\ 
\end{array} \right)
\]
is a corresponding eigenvector with $ \varphi_{i}\neq 0, \forall i=1,...,n$. 

One might think that the Jordan vectors could provide additional solutions of the associated system of differential equations in the sense of (\ref{eq118}). In the following example, we will see that this is not the case when the associated function $Q(z)$ is not a matrix polynomial. 
\hfill $\square$
\end{rem}
\begin{exm}\label{example414} Solve the system of differential equations
\begin{equation}
\label{eq428}
\begin{array}{*{20}c}
\frac{1}{u_{1}^{'}}-\frac{2}{u_{2}^{''}}=0,\\
\frac{1}{u_{1}}-\frac{2}{u_{1}^{'}}+\frac{1}{u_{1}^{''}}+\frac{1}{u_{2}}-\frac{2}{u_{2}^{'}}+\frac{1}{u_{2}^{''}}=0\\
\end{array}.
\end{equation}
\end{exm}
We seek a solution of the form (\ref{eq426}), which leads to the associated function: 
\[
Q\left( z \right)=\left( {\begin{array}{*{20}c}
\frac{1}{z} & \frac{-2}{z^{2}}\\
1-\frac{2}{z}+\frac{1}{z^{2}} & 1-\frac{2}{z}+\frac{1}{z^{2}}\\
\end{array} } \right)
\quad \wedge \quad
\det Q(z)= \frac{(z-1)^{2}(z+2)}{z^{4}}.
\]
We see that $Q(z)$ has two zeros, $\alpha_{1}=1$ of order $k_{1}=2$ and $\alpha_{2}=-2$ of order $k_{2}=1$. Therefore, if we want to find the root function at $\alpha_{1}=1$, or the corresponding Jordan vectors, we will use system (\ref{eq46}), which in this case consists of the following two equations: 

\begin{equation}
\label{eq430}
\begin{array}{*{20}c}
Q(1)\boldsymbol{\varphi}_{0} = 0,\\
Q'(1)\boldsymbol{\varphi}_{0} + Q(1)\boldsymbol{\varphi}_{1} = 0.\\
\end{array}
\end{equation}
We have
\[
Q(1)\boldsymbol{\varphi}_{0} =0 \Rightarrow 
\left( {\begin{array}{*{20}c}
1 & -2 \\
0 & 0 \\
\end{array} }\right) \left( {\begin{array}{*{20}c}
\varphi^{0}_{1}\\
\varphi^{0}_{2}\\
\end{array} }\right)=\left( {\begin{array}{*{20}c}
0\\
0\\
\end{array} }\right) \Rightarrow \boldsymbol{\varphi}_{0}=c\left( {\begin{array}{*{20}c}
2\\
1\\
\end{array} }\right), \quad c\neq 0.
\]
Note that, since $\alpha_1=1$ is an eigenvalue of $Q(z)$, the matrix $Q(1)$ is singular, as expected. 

Hence, a solution (\ref{eq426}) of the system (\ref{eq428}) is given by
\begin{equation}
\label{eq432}
\textbf{u}(t)=\left( {\begin{array}{*{20}c}
u_{1}(t)\\
u_{2}(t)\\
\end{array} }\right) =\left( {\begin{array}{*{20}c}
\frac{1}{2\varphi^{0}_{2}}\\
\frac{1}{\varphi^{0}_{2}}\\
\end{array} }\right)e^{t}, \quad \varphi^{0}_{2}\in \mathbb{C}\setminus\lbrace 0 \rbrace. 
\end{equation} 
Let us solve the second equation of (\ref{eq430}).
\[
Q'(z) = 
\begin{pmatrix}
\frac{-1}{z^{2}} & \frac{4}{z^{3}}\\
\frac{2}{z^{2}}-\frac{2}{z^{3}} & \frac{2}{z^{2}}-\frac{2}{z^{3}}\\
\end{pmatrix}
\Rightarrow 
Q'(1) = 
\begin{pmatrix}
-1 & 4\\
0 & 0\\
\end{pmatrix}.
\]
The second equation of (\ref{eq430}) gives
\[
\begin{pmatrix}
-1 & 4\\
0 & 0\\
\end{pmatrix} \left( {\begin{array}{*{20}c}
2\\
1\\
\end{array} }\right)
+
\left( {\begin{array}{*{20}c}
1 & -2 \\
0 & 0 \\
\end{array} }\right) \left( {\begin{array}{*{20}c}
\varphi^{1}_{1}\\
\varphi^{1}_{2}\\
\end{array} }\right)
= 0
\Rightarrow 
\boldsymbol{\varphi}_{1} =
\begin{pmatrix}
2(\varphi^{1}_{2}-1)\\
\varphi^{1}_{2}\\
\end{pmatrix}.
\]
If we select $\varphi_{2}^{1}=2$ we obtain
\[
\boldsymbol{\varphi}_{1} =
\begin{pmatrix}
2\\
2\\
\end{pmatrix}.
\]
One might think that formula (\ref{eq118}) could provide new solutions of the system (\ref{eq428}), following the logic of the solution (\ref{eq432}). Namely:

If we take only one vector from a Jordan chain, specifically the vector $\boldsymbol{\varphi}_{0} =
\begin{pmatrix}
\varphi^{0}_{1}\\
\varphi^{0}_{2}\\
\end{pmatrix}$, then the expression (\ref{eq118}) gives $\boldsymbol{\varphi}(z) =\begin{pmatrix}
\varphi^{0}_{1}\\
\varphi^{0}_{2}\\
\end{pmatrix}e^{t}$. This leads to the solution (\ref{eq432}), with reciprocal components forming the vector of the solution. If we take two vectors of the Jordan chain, then expression (\ref{eq118}) gives
\[
\boldsymbol{\varphi}(t)=\begin{pmatrix}
2t+2\\
t+2\\
\end{pmatrix}e^{t}.
\]
The question is: Does this lead to the new solution
\[
\textbf{u}(t)=\begin{pmatrix}
u_{1}\\
u_{2}\\
\end{pmatrix}=\left( {\begin{array}{*{20}c}
\frac{1}{2t+2}\\
\frac{1}{t+2}\\
\end{array} }\right)e^{t}
\]
of system (\ref{eq428})?
However, when we substitute the derivatives $\textbf{u}'(t)$ and $\textbf{u}''(t)$ into the first equation of system (\ref{eq428}), we see that it is not satisfied.

This gives a negative answer to the hypothesis that the Jordan vectors of the function $Q(z)$ associated with the system (\ref{eq428}) yield new solutions to the system.

Finding a solution corresponding to the eigenvalue $\alpha=-2$ is left to the reader. \hfill $\square$

Let us note that the general system of differential equations of the form (\ref{eq418}) is studied in \cite[Section 4]{B6}. In that article, the eigenvalues and eigenvectors of the associated rational matrix-valued function $Q(z)$ are obtained in a different way, rather than by solving the algebraic system (\ref{eq46}).








\vspace{.5cm}
\begin{footnotesize}
	\begin{tabular}{l}
Muhamed Borogovac\\
Boston Mutual Life\\
Actuarial Department\\
120 Royall St. Canton, MA 02021, USA\\
e-mail: {muhamed.borogovac@gmail.com }\\
\end{tabular}
\end{footnotesize}
\label{LastPage}
\end{document}